# Technical Report: A Totally Asynchronous Nesterov's Accelerated Gradient Method for Convex Optimization

Ellie Pond[*,1], April Sebok[*,2], Zachary Bell[3], and Matthew Hale[1]

*Abstract*— We present a totally asynchronous algorithm for convex optimization that is based on a novel generalization of Nesterov's accelerated gradient method. This algorithm is developed for fast convergence under "total asynchrony", i.e., allowing arbitrarily long delays between agents' computations and communications without assuming any form of delay bound. These conditions may arise, for example, due to jamming by adversaries. Our framework is block-based, in the sense that each agent is only responsible for computing updates to (and communicating the values of) a small subset of the network-level decision variables. In our main result, we present bounds on the algorithm's parameters that guarantee linear convergence to an optimizer. Then, we quantify the relationship between (i) the total number of computations and communications executed by the agents and (ii) the agents' collective distance to an optimum. Numerical simulations show that this algorithm requires 28% fewer iterations than the heavy ball algorithm and 61% fewer iterations than gradient descent under total asynchrony.

## I. INTRODUCTION

Large-scale convex optimization problems are used to model complex problems in several fields, including robotics [1], [2], [3], machine learning [4], [5], and communications [6], [7]. Large systems and/or complex tasks in these applications can lead to large convex programs, and it can be desirable to parallelize computations in order to accelerate the process of finding solutions.

Parallelized algorithms use a collection of agents to solve an optimization problem by partitioning computations among them and having the agents communicate the results of their computations with others in a network. The types of parallelized execution can be succinctly classified as (i) synchronous, (ii) partially asynchronous, or (iii) totally asynchronous. In the synchronous setting, all agents compute and communicate concurrently. However, congested bandwidth or a single slow agent can make synchrony difficult to attain. An asynchronous algorithm can be used in such cases. For partially asynchronous algorithms, the agents must compute and communicate at least once in each time interval of a prescribed length [8]. However, such bounds can be violated due to factors outside agents' control, such as adversaries jamming communications. This challenge can be alleviated by using a totally asynchronous algorithm, which allows for the delays between successive computations and communications to be unbounded for all agents [9], provided that no agent permanently stops computing and communicating.

Existing algorithms that are labeled "totally asynchronous" (or that otherwise allow for potentially unbounded probabilistic delays) include variations of gradient descent [10], [11], [12], [13], a Newton-based algorithm [14], and the heavy ball algorithm [15]. In this paper, we are motivated in part by the totally asynchronous heavy ball algorithm developed in [15], which showed faster convergence than a comparable gradient descent method. In this work, we seek even faster convergence by using Nesterov's accelerated gradient method (NAG). A feature of the heavy ball algorithm is that it converges monotonically, while NAG can converge faster than heavy ball, but potentially with oscillations [16]. In the totally asynchronous setting that we consider, the faster convergence of NAG is desirable because computations may be infrequent, which makes it critical for each computation to make as much progress as possible towards a minimizer.

Therefore, in this paper, we develop a totally asynchronous NAG algorithm that attains linear convergence to minimizers for a class of optimization problems. We apply the methodology in [9] to prove its convergence under total asynchrony. In particular, [9, pg. 431] Proposition 2.1 shows that if an algorithm is an $\infty$-norm contraction mapping, then it converges when implemented in a totally asynchronous way. The NAG algorithm itself may not be such a contraction mapping, but we we show that the application of two iterations of NAG is an $\infty$-norm contraction. We use this property to establish asymptotic convergence. Then, we establish bounds on the minimum number of computations and communications required by each agent in order for the network's iterates to be within a given distance of a minimum. Finally, we show in simulation that the totally asynchronous NAG algorithm requires up to 61% fewer iterations than gradient descent and 28% fewer than heavy ball, confirming that the desired speed up is attained.

The rest of this paper is organized as follows. Section II gives problem statements and Section III presents the totally asynchronous NAG algorithm. Section IV proves that this algorithm converges linearly, and Section V gives a convergence rate in terms of each agents' computations and communications. Next, Section VI validates the accelerated convergence of the totally asynchronous NAG algorithm in simulations. Finally, Section VII concludes.

[*]These authors contributed equally to this work.
[1]EP and MH are with the School of Electrical and Computer Engineering, Georgia Institute of Technology, Atlanta, GA 30332. Emails: {epond3,matthale}@gatech.edu.
[2]AS is with the Department of Mechanical and Aerospace Engineering, University of Florida, Gainesville, FL, USA 32611. Email: a.sebok@ufl.edu.
[3]ZB is with AFRL/RW at Eglin AFB. Email: zachary.bell.10@us.af.mil.

This work was supported by AFRL under grant FA8651-23-F-A006, AFOSR under grant FA9550-19-1-0169, and ONR under grant N00014-22-1-2435.

## II. PRELIMINARIES AND PROBLEM STATEMENTS

This section lays out notation, reviews the centralized NAG algorithm, and provides problem statements.

### A. Notation

We use $\mathbb{R}$, $\mathbb{R}_+$, and $\mathbb{N}$ to denote the real numbers, the strictly positive real numbers, and the natural numbers, respectively. We use $|\cdot|$ to denote the cardinality of a set. We use the column operator for $a, b \in \mathbb{R}^n$ defined as $\text{col}(a,b) = \begin{bmatrix} a^T & b^T \end{bmatrix}^T \in \mathbb{R}^{2n}$. We use $\Pi_{\mathcal{Z}} : \mathbb{R}^n \to \mathbb{R}^n$ to denote the orthogonal projection onto the closed, convex set $\mathcal{X} \subset \mathbb{R}^n$, i.e., $\Pi_{\mathcal{X}}[w] = \arg\min_{x \in \mathcal{X}} ||x - w||_2$. The infinity norm $||\cdot||_\infty$ is defined as $||x||_\infty = \max_{i \in \mathcal{V}} |x_i|$, where $x \in \mathbb{R}^n$ and where we define $\mathcal{V} = \{1, 2, \ldots, n\}$. We also use $\nabla_i f = \frac{\partial f}{\partial x_i}$. For an ordered pair $\chi = (v, w) \in \mathbb{R}^n \times \mathbb{R}^n$, we define $||\chi||_\infty = \max_{i \in \mathcal{V}} \max\{v_i, w_i\}$. We model the communication topology between agents as an undirected graph $\mathcal{G} = (\mathcal{V}, \mathcal{E})$, where $\mathcal{V} = \{1, 2, \ldots, n\}$ is the node set and $\mathcal{E} \subseteq \mathcal{V} \times \mathcal{V}$ is the edge set that models communication links among agents. That is, $(i, j) \in \mathcal{E}$ indicates that agents $i, j \in \mathcal{V}$ communicate with each other.

### B. Centralized Nesterov's Accelerated Gradient Method

Consider an objective function $f : \mathbb{R}^n \to \mathbb{R}$ and a constraint set $\mathcal{X} \subset \mathbb{R}^n$ (we state assumptions on them in Section III). We consider problems of the form: minimize$_{x \in \mathcal{X}} f(x)$. The centralized NAG algorithm at iteration $l \in \mathbb{N}$ updates the decision variable $x(l) \in \mathbb{R}^n$ using

$$x(l+1) = \Pi_{\mathcal{X}} \Big[ x(l) - \gamma \nabla f \big( x(l) + \lambda \big( x(l) - x(l-1) \big) \big) \\ + \lambda \big( x(l) - x(l-1) \big) \Big], \quad (1)$$

where $\lambda, \gamma \in \mathbb{R}_+$. The NAG algorithm is an accelerated gradient-based algorithm used for fast convergence. The $\gamma$ term is a step size, and here it represents the pull of "gravity," while the $\lambda$ term helps avoid overshooting the minimizer and represents "friction."

### C. Problem Statements

We consider a network of $n$ agents solving $\min_{x \in \mathcal{X}} f(x)$, where $\mathcal{X} \subset \mathbb{R}^n$ is a constraint set. The communication topology of these agents is given by a graph $\mathcal{G} = (\mathcal{V}, \mathcal{E})$.

**Remark 1.** *For simplicity, in this work we partition $x$ into scalar blocks, i.e., each agent updates a single entry of $x$, though our results will easily extend to non-scalar blocks.*

To partition the computations in (1) among agents, agent $i$ will compute successive values of $x_i$. To do so, it will need to compute $\frac{\partial f}{\partial x_i}$. However, not every other component of $x$ will explicitly appear in $\frac{\partial f}{\partial x_i}$. For example, in the objective function $f(x) = \frac{1}{2}(x_1 - x_2)^2 + \frac{1}{2}(x_3 - x_4)^2$, we have $\frac{\partial f}{\partial x_2} = -(x_1 - x_2)$, which does not depend on $x_3$ or $x_4$. Here, agent 2 does not ever need to communicate with agents 3 and 4, since agent 2's computations do not depend on those agents' decision variables and vice versa.

Generalizing this example, we see that agent $i$ does not need to communicate with any agent whose decision variable does not explicitly appear in $\frac{\partial f}{\partial x_i}$. We refer to the agents that agent $i$ must communicate with as its "essential neighbors", and we denote the set of indices of its essential neighbors as $\mathcal{V}^i \subseteq \mathcal{V}$. We emphasize here that the agents' underlying graphs do not even need to be connected as long as they communicate with their essential neighbors. In the above example with $f(x) = \frac{1}{2}(x_1 - x_2)^2 + \frac{1}{2}(x_3 - x_4)^2$, we see that agents 1 and 2 must communicate with each other and agents 3 and 4 must communicate with each other, but no communication between these two pairs is required, which allows for a graph that is not connected.

Therefore, we consider objective functions $f$ of the form

$$f(x) = \sum_{i=1}^{n} f_i(x_{\mathcal{V}^i}), \quad (2)$$

where $f_i : \mathbb{R}^{|\mathcal{V}^i|+1} \to \mathbb{R}$ and where $x_{\mathcal{V}^i}$ is a vector containing agent $i$'s decision variable and all decision variables such that $j \in \mathcal{V}^i$, i.e., the decision variables of its essential neighbors.

With this formulation, we state the problems that are the focus of the remainder of the paper.

**Problem 1.** *Given an objective function $f : \mathbb{R}^n \to \mathbb{R}$ of the form in (2) and a constraint set $\mathcal{X} \subset \mathbb{R}^n$, construct a totally asynchronous optimization algorithm based on Nesterov's accelerated gradient method that solves*

$$\min_{x \in \mathcal{X}} f(x) := \sum_{i=1}^{n} f_i(x_{\mathcal{V}^i}) \quad (3)$$

*with a network of $n$ agents.*

**Problem 2.** *Show that the totally asynchronous NAG algorithm in Problem 1 converges linearly to a minimizer.*

**Problem 3.** *Given $\epsilon > 0$, determine lower bounds on the numbers of computations and communications that each agent must execute in order for the iterates of the totally asynchronous NAG algorithm to be within distance $\epsilon$ of a minimizer for Problem 1.*

## III. TOTALLY ASYNCHRONOUS NESTEROV'S ACCELERATED GRADIENT METHOD

In this section we solve Problem 1 and formulate a totally asynchronous Nesterov's accelerated gradient (NAG) algorithm. As described in the introduction, we apply the framework from [9] to establish its convergence. Specifically, we will first show that two variations of the *synchronous* NAG algorithm satisfy certain technical conditions. Then, using [9], we will show that these convergence guarantees carry over to the totally asynchronous setting. We emphasize that all developments on a synchronous NAG algorithm will contribute to showing the convergence of the totally asynchronous NAG algorithm. Mathematically, we show that the synchronous application of two iterations of the NAG algorithm is an $\infty$-norm contraction, which will help provide the required totally asynchronous convergence guarantees.

To begin, we will make the following three assumptions about the optimization problem in (3).

**Assumption 1.** *The constraint set $\mathcal{X} \subset \mathbb{R}^n$ is nonempty, convex, and compact. The set can be decomposed as $\mathcal{X} = \mathcal{X}_1 \times \mathcal{X}_2 \times \cdots \times \mathcal{X}_n$, where $\mathcal{X}_i \subset \mathbb{R}$ for each $i \in \mathcal{V}$.*

Assumption 1 enables the parallelization of a projected update law, which will enable the constant satisfaction of set constraints, even under total asynchrony.

**Assumption 2.** *The objective function $f$ is twice continuously differentiable.*

Assumption 2 is quite common, and it guarantees the existence and continuity of the gradient and the Hessian of $f$, both of which are essential in our convergence analyses.

**Assumption 3.** *The Hessian matrix, defined as $H(x) = \nabla^2 f(x) \in \mathbb{R}^{n \times n}$, is $\mu$-diagonally dominant on $\mathcal{X} \subset \mathbb{R}^n$ for some $\mu > 0$. That is, for each $i \in \mathcal{V}$, we have the bound $H_{ii}(x) \geq \mu + \sum_{j=1, j \neq i}^{n} |H_{ij}(x)|$ for all $x \in \mathcal{X}$.*

Assumption 3 is standard in the context of totally asynchronous algorithms. Intuitively, this assumption asserts that for agent $i \in \mathcal{V}$, its computations depend more on its own decision variables than on the rest of the agents' decision variables. In [9], it is noted that some form of Hessian diagonal dominance is typically needed for convergence of totally asynchronous algorithms, and we therefore use it here.

**Remark 2.** *Assumption 3 implies that $f$ is $\mu$-strongly convex. Therefore it has a unique minimizer over $\mathcal{X}$, which we denote by $x^\star = \mathrm{col}(x_1^\star, \ldots, x_n^\star)$.*

In its centralized form in (1), Nesterov's accelerated gradient method depends on the iterate at time $l$, namely $x(l)$, and the iterate at the previous time step $l-1$, which is $x(l-1)$. We write $y(l) = x(l-1)$ to denote this latter iterate. For the distributed solution of Problem 1, each agent will store a local copy of the full decision vector in its onboard memory. Onboard agent $i$ at time $l$, this decision variable is denoted $z^i(l) = (x^i(l), y^i(l)) \in \mathcal{Z}$, where $\mathcal{Z} := \mathcal{X} \times \mathcal{X}$ and where the superscript on each vector denotes that it is stored onboard agent $i$.

Over time, agent $i$ computes updates to $z_i^i(l) = (x_i^i(l), y_i^i(l)) \in \mathcal{Z}_i$, where we define $\mathcal{Z}_i := \mathcal{X}_i \times \mathcal{X}_i$. The subscripts indicate that the terms $z_i^i$, $x_i^i$, and $y_i^i$ are agent $i$'s local copies of its own decision variable. Using this notation, at time $l$ agent $i$'s local copy of $z$ is denoted

$$z^i(l) = (x^i(l), y^i(l)) = \left(\mathrm{col}(x_1^i(l), \ldots, x_i^i(l), \ldots, x_n^i(l)), \mathrm{col}(y_1^i(l), \ldots, y_i^i(l), \ldots, y_n^i(l))\right).$$

For any agent $m \notin \mathcal{V}^i$, agent $i$ can set $x_m^i$ and $y_m^i$ to any values over time, since these values do not affect its computations and will not be changed by any communications.

### A. The Single-Step Synchronous Method

In this sub-section, we establish the first of two variations of the synchronous NAG algorithm that we use in Section III-C for the analysis of the totally asynchronous NAG algorithm. These algorithms have simultaneous computations and simultaneous communications among essential neighbors.

We will refer to the following update law for $z_i^i(l) = (x_i^i(l), y_i^i(l))$ as the "single-step synchronous update law". This update law is defined by $\tilde{u}_x^i$ and $\tilde{u}_y^i$, given by

$$x_i^i(l+1) = \tilde{u}_x^i(x^i(l), y^i(l)) \quad (4)$$
$$= \Pi_{\mathcal{X}_i}\left[x_i^i(l) - \gamma \nabla_i f\left(x^i(l) + \lambda(x^i(l) - y^i(l))\right) + \lambda(x_i^i(l) - y_i^i(l))\right]$$

$$y_i^i(l+1) = \tilde{u}_y^i(x^i(l), y^i(l)) \quad (5)$$
$$= x_i^i(l)$$

for all $l \in \mathbb{N}$ and $i \in \mathcal{V}$. In this update law, we see that at time $l \in \mathbb{N}$, one iteration of NAG is performed and stored in the $x_i^i(l+1)$ variable, while the $y_i^i(l+1)$ variable stores the value of $x_i^i(l)$. For simplicity of notation, let $\tilde{u}^i : \mathcal{Z} \to \mathcal{Z}_i$ denote the single-step synchronous update

$$\left(x_i^i(l+1), y_i^i(l+1)\right) = \tilde{u}^i(x^i(l), y^i(l)) \quad (6)$$
$$= \left(\tilde{u}_x^i(x^i(l), y^i(l)), \tilde{u}_y^i(x^i(l), y^i(l))\right).$$

We will also denote (6) by $z_i^i(l+1) = \tilde{u}^i(z^i(l))$ for conciseness. The following lemma establishes that this update law has a single fixed point, which is the solution to Problem 1.

**Lemma 1.** *Consider Problem 1, and let $\mathcal{X}$ satisfy Assumption 1, $f$ satisfy Assumption 2, and the Hessian $H(x) = \nabla^2 f(x)$ satisfy Assumption 3 with $\mu > 0$. Let $x^\star$ be as defined in Remark 2. Define the points $z^\star = (x^\star, x^\star) \in \mathcal{Z}$ and $z_i^\star = (x_i^\star, x_i^\star) \in \mathcal{Z}_i$. Then, the point $z^\star$ is a fixed point of the single-step synchronous update law in (6), in the sense that $z_i^\star = \tilde{u}^i(z^\star)$ for all $i \in \mathcal{V}$.*

*Proof.* See Appendix A. □

In the synchronous algorithm, all agents update and communicate at each time step $l \in \mathbb{N}$. Agent $i$ updates its decision variables according to (6) and communicates this update with agents in its essential neighborhood, indexed by $j \in \mathcal{V}^i$. That is, at time step $l \in \mathbb{N}$ agent $i$ sends $(x_i^i(l), y_i^i(l))$ to all agents $j \in \mathcal{V}^i$. All agents $j \in \mathcal{V}^i$ incorporate this communication into their own local state vector by setting $(x_i^j(l), y_i^j(l)) \leftarrow (x_i^i(l), y_i^i(l))$. For the agents $m \notin \mathcal{V}^i$, the entries in their local state vector remain the same as at the previous time step, i.e., $(x_i^m(l), y_i^m(l)) \leftarrow (x_i^m(l-1), y_i^m(l-1))$.

We define the true state of the network at time $l \in \mathbb{N}$ to be the vector $z^{\mathrm{true}}(l) = (x^{\mathrm{true}}(l), y^{\mathrm{true}}(l))$, where $x^{\mathrm{true}}(l) = \mathrm{col}(x_1^1(l), x_2^2(l), \ldots, x_n^n(l))$ and $y^{\mathrm{true}}(l) = \mathrm{col}(y_1^1(l), y_2^2(l), \ldots, y_n^n(l))$. These vectors contain each agent's latest value of their own decision variable. To establish the convergence of the single-step

synchronous NAG algorithm, we will use the map

$$\hat{u}_{\text{true}}^i(z^{\text{true}}(l)) = x_i^i(l) - \gamma \nabla_i f\big(x(l) + \lambda(x(l) - y(l))\big) \\ + \lambda(x_i^i(l) - y_i^i(l)), \quad (7)$$

which models the changes in the true state of the network at each iteration.

**Theorem 1.** *Consider Problem 1, and let $\mathcal{X}$ satisfy Assumption 1, $f$ satisfy Assumption 2, and the Hessian $H(x) = \nabla^2 f(x)$ satisfy Assumption 3 with $\mu > 0$. For each $\gamma \in \left(0, \frac{1}{\max_{i \in \mathcal{V}} \max_{\eta \in \mathcal{X}} |H_{ii}(\eta)|}\right)$ and $\lambda \in \left(0, \frac{\gamma \mu}{2(1-\gamma\mu)}\right)$, the iterates of the synchronous single-step algorithm from the initial state $z(0) \in \mathcal{Z}$ satisfy $\|z(l+1) - z^\star\|_\infty \leq \alpha \|z(l-1) - z^\star\|_\infty$ for all $l \in \mathbb{N}$, where $\alpha = \max\{\alpha_1, \alpha_2\}$ and*

$$\alpha_1 = \big(1 + \lambda - \gamma\mu(1+\lambda)\big)^2 + \lambda(1-\gamma\mu) \quad (8) \\ + \lambda(1-\gamma\mu)\big(1 + \lambda - \gamma\mu(1+\lambda)\big)$$

$$\alpha_2 = 1 - \gamma\mu + 2\lambda(1-\gamma\mu), \quad (9)$$

*with $\alpha_1, \alpha_2 \in [0, 1)$.*

*Proof.* See Appendix B. □

**Remark 3.** *Theorem 1 proves that the synchronous single-step NAG update law is contractive with respect to the $\infty$-norm over two time steps, i.e., from time $l-1$ to time $l+1$.*

### B. The Double-Step Synchronous Method

In this section, we continue with the procedure for proving totally asynchronous convergence outlined in [9]. This process requires proving that the synchronous variation of the algorithm is contractive with respect to the infinity-norm, which we have shown is the case over two time steps in the previous subsection. For this reason, we now define what we term the "double-step synchronous update law", which performs two steps of NAG per iteration and hence is contractive at every iteration. As well, we will prove an additional three-part lemma in regard to the double-step NAG law that is required for totally asynchronous convergence.

We will now use the variable $k \in \mathbb{N}$ to represent discretized time rather than $l \in \mathbb{N}$. In the $k$ time-scale, the step from $k \to k+1$ is equivalent to $l \to l+2$ in the $l$ time-scale. We make this change in order to make it clear that the NAG algorithm in (1) is applied twice every time an agent performs a computation. For each $i \in \mathcal{V}$, agent $i$ updates in the double-step NAG update law take the form

$$x_i^i(k+1) = u_x^i(x^i(k), y^i(k+1)) \quad (10) \\ = \Pi_{\mathcal{X}_i} \Big[ y_i^i(k+1) + \lambda\big(y_i^i(k+1) - x_i^i(k)\big) \\ - \gamma \nabla_i f\big(y^i(k+1) + \lambda(y^i(k+1) - x^i(k))\big) \Big]$$

$$y_i^i(k+1) = u_y^i(x^i(k), y^i(k)) \quad (11) \\ = \Pi_{\mathcal{X}_i} \Big[ x_i^i(k) + \lambda\big(x_i^i(k) - y_i^i(k)\big) \\ - \gamma \nabla_i f\big(x^i(k) + \lambda(x^i(k) - y^i(k))\big) \Big]$$

for all $k \in \mathbb{N}$. The local state vector $y^i(k+1)$ that is used in (10) is defined as $y^i(k+1) = \text{col}(y_1^i(k), y_2^i(k), \ldots, y_i^i(k+1), \ldots, y_n^i(k))$, where the newly updated $y_i^i$ is stored in the $i^{\text{th}}$ location and all other entries remain the same as they were at time $k$, i.e., $y_j^i(k+1) = y_j^i(k)$ for $j \neq i$.

**Remark 4.** *In (10), the variable $y_i^i(k+1)$ is referenced. Though this may appear non-recursive, the update for $x_i^i$ can be expressed completely in terms of $(x^i(k), y^i(k))$ using the update law for $y_i^i(k+1)$ in (11). However, we choose to refer to the $(k+1)^{\text{th}}$ iteration of $y_i^i$ for brevity of expression. Nonetheless, both (10) and (11) can be computed simultaneously over one time step $k \in \mathbb{N}$.*

Let $u^i : \mathcal{Z} \to \mathcal{Z}_i$ be the double-step synchronous NAG update for agent $i$, defined as

$$(x_i^i(k+1), y_i^i(k+1)) = u^i(x^i(k), y^i(k)) \\ = \big(u_x^i(x^i(k), y^i(k+1)), u_y^i(x^i(k), y^i(k))\big). \quad (12)$$

With an abuse of notation, we will alternately write (12) as $z_i^i(k+1) = u^i(z^i(k))$ for conciseness. As in Section III-A, every agent updates and communicates with their essential neighbors at every time step $k \in \mathbb{N}$ for the synchronous double-step algorithm. This update law has the same fixed point at the synchronous single-step algorithm.

**Lemma 2.** *Consider Problem 1, and let $\mathcal{X}$ satisfy Assumption 1, $f$ satisfy Assumption 2, and the Hessian $H(x) = \nabla^2 f(x)$ satisfy Assumption 3. Then, for each $i \in \mathcal{V}$, the minimizer $z^\star \in \mathcal{Z}$ is a fixed point of the double-step synchronous update in (12), in the sense that $z_i^\star = u^i(z^\star)$ for all $i \in \mathcal{V}$.*

*Proof.* This proof follows the procedure given for the proof of Lemma 1 and is therefore omitted for brevity. □

For ease of notation, we define $h : \mathcal{Z} \to \mathcal{Z}$ to be the map

$$h(z) := \text{col}\Big(\big(u_x^1(z), u_x^2(z), \ldots, u_x^n(z)\big)^T, \\ \big(u_y^1(z), u_y^2(z), \ldots, u_y^n(z)\big)^T\Big), \quad (13)$$

which is equivalent to one iterate of the double-step synchronous NAG update law. The point $z^\star \in \mathcal{Z}$ is a fixed point of $h$ due to Lemma 2.

The following lemma defines and analyzes a collection of sets $\{Z(k)\}_{k \in \mathbb{N}}$ that will be used in the analysis of the forthcoming totally asynchronous NAG algorithm. These sets will essentially serve as Lyapunov sub-level sets in that analysis, and this collection of sets must satisfy three conditions: (i) the Lyapunov-Like Condition (LLC), (ii) the Synchronous Convergence Condition (SCC), and (iii) the Box Condition (BC), which we define next.

**Lemma 3.** *Consider Problem 1, and let $\mathcal{X}$ satisfy Assumption 1, $f$ satisfy Assumption 2, and the Hessian $H(x) = \nabla^2 f(x)$ satisfy Assumption 3 with $\mu > 0$. Let $z(0) \in \mathcal{Z}$ be*

given, and define the set $\mathcal{Z}(k)$ as

$$\mathcal{Z}(k) = \{v \in \mathcal{Z} : ||v - z^\star||_\infty \leq \alpha^k ||z(0) - z^\star||_\infty\},$$

where $\alpha = \max\{\alpha_1, \alpha_2\}$ is from Theorem 1. Then, for every $k \in \mathbb{N}$, the set $\mathcal{Z}(k)$ satisfies the following three properties:

1) (LLC) The set containment rule

$$\cdots \subset \mathcal{Z}(k+1) \subset \mathcal{Z}(k) \subset \cdots \subset \mathcal{Z}(0) = Z \quad (14)$$

holds for all $k \in \mathbb{N}$.

2) (SCC) For the mapping $h$ in (13), given a point $z \in \mathcal{Z}(k)$, we have $h(z) \in \mathcal{Z}(k+1)$ for all $k \in \mathbb{N}$. As well, if $\{z_k\}_{k \in \mathbb{N}}$ is a sequence such that each $z_k \in \mathcal{Z}(k)$ for each $k \in \mathbb{N}$, then $\lim_{k \to \infty} z_k = z^\star$, where $z^\star = (x^\star, x^\star)$ is the fixed point of $h$.

3) (BC) For all $k \in \mathbb{N}$ and $i \in \mathcal{V}$, there are sets $\mathcal{Z}_i(k) \subset \mathcal{Z}_i$ such that $\mathcal{Z}(k) = \mathcal{Z}_1(k) \times \mathcal{Z}_2(k) \times \cdots \times \mathcal{Z}_n(k)$.

*Proof.* We prove items 1) through 3) in order.

1) For an arbitrary $k \in \mathbb{N}$ we consider the sets

$$\mathcal{Z}(k) = \{v \in \mathcal{Z} : ||v - z^\star||_\infty \leq \alpha^k ||z(0) - z^\star||_\infty\}$$

and

$$\mathcal{Z}(k+1) = \{v \in \mathcal{Z} : ||v - z^\star||_\infty \leq \alpha^{k+1} ||z(0) - z^\star||_\infty\}.$$

Since $\alpha < 1$ due to Theorem 1, we see that $\alpha^{k+1} \leq \alpha^k$. Then $v \in \mathcal{Z}(k+1)$ satisfies

$$||v - z^\star||_\infty \leq \alpha^{k+1} ||z(0) - z^\star||_\infty \leq \alpha^k ||z(0) - z^\star||_\infty.$$

Then $v \in \mathcal{Z}(k)$ and $\mathcal{Z}(k+1) \subseteq \mathcal{Z}(k)$. Because $k \in \mathbb{N}$ was abitrary, we see that $\cdots \subset \mathcal{Z}(k+2) \subset \mathcal{Z}(k+1) \subset \mathcal{Z}(k) \subset \cdots Z$ for all $k \in \mathbb{N}$.

2) (SCC) By definition of $\mathcal{Z}(k)$, having $z \in \mathcal{Z}(k)$ implies that $||z - z^\star||_\infty \leq \alpha^k ||z(0) - z^\star||_\infty$. The mapping $h$ implements two steps of the NAG algorithm, and Theorem 1 shows that two steps of the NAG algorithm result in a contraction of the $\infty$-norm distance to the minimizer by a factor of $\alpha$. That is, Theorem 1 implies that

$$||h(z) - z^\star||_\infty \leq \alpha ||z - z^\star||_\infty \leq \alpha^{k+1} ||z(0) - z^\star||_\infty.$$

Then $h(z) \in \mathcal{Z}(k+1)$, as desired.

Now consider a sequence $\{z_k\}_{k \in \mathbb{N}}$, where $z_k \in \mathcal{Z}(k)$ for each $k \in \mathbb{N}$. From the LLC, we know that (14) holds and $\mathcal{Z}(k) \subset \mathcal{Z}(k-1)$ for all $k \in \mathbb{N}$. By inspection each set $\mathcal{Z}(k)$ is closed. Therefore, we have, $\lim_{k \to \infty} \mathcal{Z}(k) = \bigcap_{k \in \mathbb{N}} \mathcal{Z}(k)$ by [17]. Thus, here we find

$$\lim_{k \to \infty} z_k \in \lim_{k \to \infty} Z(k) = \bigcap_{k \in \mathbb{N}} \mathcal{Z}(k) = \{z^\star\}.$$

As well, in Lemma 2 it was shown that $z^\star \in \mathcal{Z}$ furnishes a fixed point of each $u^i$, i.e., $z_i^\star = u^i(z^\star)$ for all $i \in \mathcal{V}$. Using the definition of $h$ in (13), we see that $\lim_{k \to \infty} z_k = z^\star$, and $z^\star$ is indeed a fixed point of $h$.

3) (BC) By the definition of the infinity norm, for $v \in \mathcal{Z}(k)$ we have $||v - z^\star||_\infty = \max_{i \in \mathcal{V}} |v_i - z_i^\star|$. By definition of $\mathcal{Z}(k)$, this equality then implies that $|v_i - z_i^\star| \leq \alpha^k ||z(0) - z^\star||_\infty$ for all $i \in \mathcal{V}$. Then we see that, for all $i \in \mathcal{V}$, we have $v_i \in \mathcal{Z}_i(k)$, where $\mathcal{Z}_i(k) = \{v_i \in \mathcal{Z}_i : |v_i - z_i^\star| \leq \alpha^k ||z(0) - z^\star||_\infty\}$, and thus $\mathcal{Z}(k) = \mathcal{Z}_1(k) \times \mathcal{Z}_2(k) \times \cdots \times \mathcal{Z}_n(k)$. $\square$

## C. Totally Asynchronous NAG Algorithm

In this section, we build upon the previous subsections and develop the totally asynchronous NAG algorithm. Under total asynchrony, at any particular timestep each agent may or may not perform a computation. To keep track of when each agent does so, we let $K^i \subseteq \mathbb{N}$ be the set of time steps at which agent $i \in \mathcal{V}$ updates its own decision variables. We emphasize that the set $K^i$ for each $i$ is not known to the agents and is simply used to facilitate analysis.

When performing a computation, agent $i$ uses the same update law in the totally asynchronous setting as in the synchronous double-step setting, namely (12), but it can do so with any timing and without coordinating that timing with any other agent. To faithfully implement that update law, any communications that agent $i$ receives from another agent at time $k$ are not incorporated into the local state vector $y^i(k+1)$ until after agent $i$ has computed an update to both $x_i^i$ and $y_i^i$. Therefore, at each time $k \in K^i$, agent $i$ updates its decision variables with $\bigl(x_i^i(k+1), y_i^i(k+1)\bigr) \leftarrow u^i\bigl(x^i(k), y^i(k)\bigr)$.

After an update is computed at a time $k \in K^i$, agent $i$ sends its updated decision variables $\bigl(x_i^i(k+1), y_i^i(k+1)\bigr)$ to each agent $j \in \mathcal{V}^i$. However, in the totally asynchronous setting, agent $i$ may not communicate every one of its updated decision variable values over time. Moreover, there can be communication delays between agent $i$ sending an updated value of its iterates and agent $j$ receiving it. To model these communications, we let $R_j^i \subseteq \mathbb{N}$ contain the set of times at which agent $i$ receives a communication from its essential neighbor $j$; if $j$ is not an essential neighbor of agent $i$, then they do not communicate and we have $R_j^i = \emptyset$. As with $K^i$, we emphasize that the sets $R_j^i$ are not known to the agents and are only used to facilitate analysis.

At any time $k \in R_j^i$, agent $i$ uses the values it receives from agent $j$ to overwrite the previous values of agent $j$'s iterates that it had received from agent $j$. Formally, agent $i$ executes the operation $\bigl(x_j^i(k), y_j^i(k)\bigr) \leftarrow \bigl(x_j^j\bigl(\tau_j^i(k)\bigr), y_j^j\bigl(\tau_j^i(k)\bigr)\bigr) \in \mathcal{Z}_j$. Here, the notation $\tau_j^i(k) \in K^j$ denotes the time at which agent $j$ originally computed the values of $x_j^j$ and $y_j^j$ that agent $i$ has onboard at time $k$. That is, at all times $k \in \mathbb{N}$, we define $\tau_j^i$ to be the earliest time in $K^j$ so that

$$\bigl(x_j^i(k), y_j^i(k)\bigr) = \bigl(x_j^j\bigl(\tau_j^i(k)\bigr), y_j^j\bigl(\tau_j^i(k)\bigr)\bigr)$$

holds. For agents $m \notin \mathcal{V}^i$, the entries in the local state vector $\bigl(x_m^i(k+1), y_m^i(k+1)\bigr)$ remain constant since agents $i$ and $m$ do not communicate and the values of agent $m$'s decision variables do not affect agent $i$'s computations.

The full NAG algorithm with totally asynchronous computations and communications is shown in Algorithm 1, and this algorithm solves Problem 1.

**Algorithm 1** Totally Asynchronous NAG Algorithm

**Input:** For $i \in \mathcal{V}$ select an arbitrary initial state $z^i(0) \in \mathcal{Z}$.

1 **for** $k \in \mathbb{N}$ **do**
2     **for** $i \in \mathcal{V}$ **do**
3         **if** $k \in K^i$ **then**
4             $\left(x_i^i(k+1), y_i^i(k+1)\right) \leftarrow u^i\left(x^i(k), y^i(k)\right)$
5             **if** $j \in \mathcal{V}^i$ **then**
6                 Send $\left(x_i^i(k+1), y_i^i(k+1)\right)$ to agent $j$
7             **end**
8         **end**
9         **for** $j \in \mathcal{V}^i$ **do**
10             **if** $k \in R_j^i$ **then**
11                 $\left(x_j^i(k), y_j^i(k)\right) \leftarrow \left(x_j^j(\tau_j^i(k)), y_j^j(\tau_j^i(k))\right)$
12             **end**
13         **end**
14         **if** $m \notin \mathcal{V}^i$ **then**
15             $\left(x_m^i(k), y_m^i(k)\right) \leftarrow \left(x_m^i(k-1), y_m^i(k-1)\right)$
16         **end**
17     **end**
18 **end**

## IV. CONVERGENCE RATE

In this section we will show that the totally asynchronous NAG Algorithm 1 converges linearly, which will solve Problem 2. To begin, we state an assumption regarding agents' computations and communications and introduce the notion of an "operation cycle".

**Assumption 4** ([9]). *For each agent $i \in \mathcal{V}$ and each essential neighbor $j \in \mathcal{V}^i$, the update set $K^i$ and the communication set $R_j^i$ are infinite. Moreover, if $\{k_s\}_{s \in \mathbb{N}}$ is an increasing sequence of times in $K^i$, then $\lim_{s \to \infty} \tau_j^i(k_s) = \infty$ for all $j \in \mathcal{V}$.*

This assumption ensures that no agent will ever cease computing or communicating indefinitely, though it allows for the delays between successive computations and communications to be arbitrarily long.

For the totally asynchronous setting, we will define an "operation cycle" to analyze Algorithm 1's convergence. Over time, such cycles are tracked by a counter, and we say that ops($k$) cycles have been completed by time $k$. Initially, when $k = 0$, the operation cycle counter ops(0) = 0. Suppose some time $k' \in \mathbb{N}$ is the first point in time by which (i) every agent $i \in \mathcal{V}$ has computed an update to their decision variables, (ii) they have sent these updated values to all of their essential neighbors, and (iii) every essential neighbor $j \in \mathcal{V}^i$ has received these values and incorporated them into their local state vector. Then, at that time, we have ops($k'$) = 1.

Incrementing this value reflects the fact that one operation cycle has been completed once (i)-(iii) have been completed.

After timestep $k'$, the next communication cycle is completed after (i)-(iii) have been completed again. Suppose this happens at time $\ell \in \mathbb{N}$. Then we have ops($k$) = 1 for all $k' \leq k < \ell$, and ops($\ell$) = 2. The value of ops($\cdot$) will remain equal to 2 until the next cycle is completed.

Note that one or more agents may compute and communicate more than one time per cycle. However, an operation cycle is not completed until *every* agent has performed a computation, communicated it to its essential neighbors, and received a communication from each of its essential neighbors at least once. At a time where ops($k$) = 0, each agent's local state satisfies $\left(x^i(k), y^i(k)\right) \in \mathcal{Z}(0) = \mathcal{Z}$. Once the first operation cycle is completed at time $k'$, we have ops($k'$) = 1, and based on that definition every agent's local state satisfies $\left(x^i(k'), y^i(k')\right) \in \mathcal{Z}(1)$.

We now establish an invariance result that is needed to derive the convergence rate of the algorithm.

**Lemma 4.** *Consider Problem 1, and let $\mathcal{X}$ satisfy Assumption 1, $f$ satisfy Assumption 2, and the Hessian $H(x) = \nabla^2 f(x)$ satisfy Assumption 3. Given an initial state $z^i(0) \in \mathcal{Z}$ for every $i \in \mathcal{V}$, the set $\mathcal{Z}(k)$ is forward invariant for the totally asynchronous NAG algorithm in Algorithm 1. That is, for all $i \in \mathcal{V}$, once $z^i(l) \in \mathcal{Z}(k)$ for some $l \in \mathbb{N}$, it holds that $z^i(p) \in \mathcal{Z}(k)$ for all $p \geq l$.*

*Proof.* See Appendix C. □

**Lemma 5.** *Consider Problem 1, and let $\mathcal{X}$ satisfy Assumption 1, $f$ satisfy Assumption 2, the Hessian $H(x) = \nabla^2 f(x)$ satisfy Assumption 3, and let Assumption 4 hold. Then, for each $i \in \mathcal{V}$, the minimizer $z^\star \in \mathcal{Z}$ is a fixed point of the totally asynchronous NAG update law (12), in the sense that $z_i^\star = u^i(z^\star)$ for all $i \in \mathcal{V}$.*

*Proof.* This proof follows the procedure given for the proof of Lemma 1 and is therefore omitted for brevity. □

Now we present the theorem that establishes that the totally asynchronous NAG algorithm converges linearly in the value of ops($k$), which solves Problem 2.

**Theorem 2.** *Consider Problem 1, and let $\mathcal{X}$ satisfy Assumption 1, $f$ satisfy Assumption 2, the Hessian $H(x) = \nabla^2 f(x)$ satisfy Assumption 3 with $\mu > 0$, and let Assumption 4 hold. For each $\gamma \in \left(0, \frac{1}{\max_{i \in \mathcal{V}} \max_{\eta \in \mathcal{X}} |H_{ii}(\eta)|}\right)$ and $\lambda \in \left(0, \frac{\gamma \mu}{2(1 - \gamma \mu)}\right)$, the totally asynchronous NAG algorithm in Algorithm 1 in which agent $i$ is initialized with the initial condition $z^i(0) \in \mathcal{Z}$ for all $i \in \mathcal{V}$ satisfies*

$$\max_{i \in \mathcal{V}} \left|\left|z^i(k) - z^\star\right|\right|_\infty \leq \alpha^{\text{ops}(k)} \max_{i \in \mathcal{V}} \left|\left|z^i(0) - z^\star\right|\right|_\infty$$

*for all $k \in \mathbb{N}$, where $\alpha = \max\{\alpha_1, \alpha_2\}$ with $\alpha_1 \in [0, 1)$ defined in (8) and $\alpha_2 \in [0, 1)$ defined in (9).*

*Proof.* Recall that each agent $i \in \mathcal{V}$ must update its decision variables and communicate them to all agents in $j \in \mathcal{V}^i$ at least once during an operation cycle. Consider an initial state $z^i(0) \in \mathcal{Z}(0) = \mathcal{Z}$. For each $i \in \mathcal{V}$, suppose that agent $i$ updates its decision variables for the first time at some $k_i \in$

$K^i$. Then $\big(x_i^i(k_i+1), y_i^i(k_i+1)\big) \leftarrow u^i\big(x^i(k_i), y^i(k_i)\big)$. By the SCC and the BC in Lemma 3, $\big(x_i^i(k_i+1), y_i^i(k_1+1)\big) \in \mathcal{Z}_i(1)$. At time $k' = \max_{i \in \mathcal{V}} k_i + 1$, it is then true that $\big(x_i^i(k'), y_i^i(k')\big) \in \mathcal{Z}_i(1)$ holds for all agents $i \in \mathcal{V}$ by Assumption 4. However, we still have $\big(x^i(k'), y^i(k')\big) \in \mathcal{Z}(0)$ at time $k'$, because $\big(x_i^i(k'), y_i^i(k')\big) \in Z_i(1) \subset \mathcal{Z}_i(0)$ by the LLC in Lemma 3 and $\big(x_j^i(k'), x_j^i(k')\big) \in \mathcal{Z}_j(0)$ prior to agent $i$ receiving communications from agent $j$. This statement holds for both $j \in \mathcal{V}^i$ and $j \notin \mathcal{V}^i$.

After each agent $i \in \mathcal{V}$ has completed an update, it sends this information to agents $j \in \mathcal{V}^i$ after some amount of time due to Assumption 4. Suppose that at time $k_{i,j} \in R_j^i$, agent $i \in \mathcal{V}$ receives $\big(x_j^j\big(\tau_j^i(k_{i,j})\big), y_j^j\big(\tau_j^i(k_{i,j})\big)\big)$ from essential neighbor $j \in \mathcal{V}^i$ for the first time. Then, there exists a time $k'' \geq k'$ defined as $k'' = \max_{i \in \mathcal{V}} \max_{j \in \mathcal{V}^i} k_{i,j}$ by which every agent has received a communication from each of its essential neighbors. That is, $\big(x_j^i(k''), y_j^i(k'')\big) \in \mathcal{Z}_j(1)$ for all $i \in \mathcal{V}$ and all $j \in \mathcal{V}^i$. Therefore, the local vectors $\big(x^i(k''), y^i(k'')\big) \in \mathcal{Z}(1)$ for all $i \in \mathcal{V}$. With this, the first operation cycle has been completed at time $k''$ and we have $\mathrm{ops}(k'') = 1$.

By definition of $\mathcal{Z}(1)$ in Lemma 3, the local vectors $z^i(k'') = \big(x^i(k''), y^i(k'')\big)$ satisfy
$$\big|\big|z^i(k'') - z^\star\big|\big|_\infty \leq \alpha \big|\big|z^i(0) - z^\star\big|\big|_\infty$$
for all $i \in \mathcal{V}$, where $\alpha^{\mathrm{ops}(k'')} = \alpha$ here since $\mathrm{ops}(k'') = 1$. By induction, we may generalize this relationship to times $k$ beyond $k''$ to find that
$$\big|\big|z^i(k) - z^\star\big|\big|_\infty \leq \alpha^{\mathrm{ops}(k)} \big|\big|z^i(0) - z^\star\big|\big|_\infty$$
for all $i \in \mathcal{V}$.

Thus, over the entire network, linear convergence holds, i.e.,
$$\max_{i \in \mathcal{V}} \big|\big|z^i(k) - z^\star\big|\big|_\infty \leq \alpha^{\mathrm{ops}(k)} \max_{i \in \mathcal{V}} \big|\big|z^i(0) - z^\star\big|\big|_\infty$$
for all $k \in \mathbb{N}$. $\square$

## V. OPERATION COMPLEXITY

In this section, we will leverage the convergence rate established in Theorem 2 and the network properties of $\mathcal{G} = (\mathcal{V}, \mathcal{E})$ to quantify the operation complexity of agents' convergence. As defined previously, a single operation cycle occurs when all agents update their decision variables, communicate those values to their essential neighbors, and receive updated values of their essential neighbors' decision variables. Rephrased, an operation cycle consists of at least $|\mathcal{V}|$ computation events (since every agent must perform at least one computation per operation cycle) and at least $2|\mathcal{E}|$ communication events (since every agent must send information to and receive information from each of their essential neighbors). On the agent level, agent $i$ performs a computation at least one time, sends information to essential neighbors at least $|\mathcal{V}^i|$ total times, and receives information from essential neighbors at least $|\mathcal{V}^i|$ total times.

Our goal is to establish bounds on the number of computations and communications that are required in order for the network's solution to be within an $\infty$-norm ball of radius $\epsilon$ that is centered on the minimizer of (3). That is, we will bound the number of operations required for agent $i$ to reach a point $z^i(k)$ such that $\big|\big|z^i(k) - z^\star\big|\big|_\infty \leq \epsilon$ for all $i$. Of course, the minimzer $z^\star$ is unknown in general, and for this reason we use $D_0 := \max_{v_1, v_2 \in \mathcal{Z}} ||v_1 - v_2||_\infty$ as an upper bound on the initial distance of agents' iterates to a minimizer.

**Theorem 3.** *Consider Problem 1, and let $\mathcal{X}$ satisfy Assumption 1, $f$ satisfy Assumption 2, the Hessian $H(x) = \nabla^2 f(x)$ satisfy Assumption 3 with $\mu > 0$, and let Assumption 4 hold. Given $\gamma \in \left(0, \frac{1}{\max_{i \in \mathcal{V}} \max_{\eta \in \mathcal{X}} |H_{ii}(\eta)|}\right)$ and $\lambda \in \left(0, \frac{\gamma \mu}{2(1-\gamma\mu)}\right)$, in order for agent $i$ to be within $\epsilon$ of the minimizer $z^\star$ for all $i \in \mathcal{V}$, i.e., to attain*
$$\max_{i \in \mathcal{V}} \big|\big|z^i(k) - z^\star\big|\big|_\infty \leq \epsilon, \tag{15}$$
*agent $i$ must have performed at least $\beta$ computations and communicated at least $\beta |\mathcal{V}^i|$ times, where $\beta = \frac{\log(\epsilon/D_0)}{\log(\alpha)}$.*

*Proof.* From Theorem 2, the network as a whole satisfies
$$\max_{i \in \mathcal{V}} \big|\big|z^i(k) - z^\star\big|\big|_\infty \leq \alpha^{\mathrm{ops}(k)} \max_{i \in \mathcal{V}} \big|\big|z^i(0) - z^\star\big|\big|_\infty$$
for all $k \in \mathbb{N}$. Then,
$$\max_{i \in \mathcal{V}} \big|\big|z^i(k) - z^\star\big|\big|_\infty \leq \alpha^{\mathrm{ops}(k)} \max_{v_1, v_2 \in \mathcal{Z}} ||v_1 - v_2||_\infty$$
$$= \alpha^{\mathrm{ops}(k)} D_0.$$

We therefore seek to enforce $\alpha^{\mathrm{ops}(k)} D_0 \leq \epsilon$.

Solving for $\mathrm{ops}(k)$ results in $\mathrm{ops}(k) \geq \frac{\log(\epsilon/D_0)}{\log(\alpha)}$. We define $\beta := \frac{\log(\epsilon/D_0)}{\log(\alpha)}$. In every operation cycle, each agent updates at least once and communicates at least $|\mathcal{V}^i|$ times. Thus, in order for the network to satisfy (15), each agent must have completed at least $\beta$ updates and $\beta |\mathcal{V}^i|$ communications during the execution of Algorithm 1. $\square$

## VI. SIMULATIONS

To demonstrate the effectiveness of Algorithm 1, we provide a simulation in MATLAB. We consider a network of 10 agents with the objective function
$$f(x) = \frac{3}{10} \sum_{i=1}^n (x_i^i)^2 + \frac{1}{200} \sum_{i=1}^n \sum_{\substack{j=1 \\ j \neq i}}^n (x_i^i - x_j^i)^2,$$
where $\mu = 0.3$. Our constraint set for each agent $i \in \mathcal{V}$ is $\mathcal{Z}_i = [1, 10]$. We define the hyperparameters $\lambda$ and $\gamma$ using $\lambda = \frac{\gamma \mu}{2(1-\gamma\mu)} = 0.058$ and $\gamma = \frac{1}{\max_{i \in \mathcal{V}} \max_{\eta \in \mathcal{X}} |H_{ii}(\eta)|} = 0.345$. The initial conditions for each agent are $z_0^i = (10 \cdot \mathbf{1}, 10 \cdot \mathbf{1})$, where $\mathbf{1}$ is the 10-dimensional vector of ones. We use a uniform probability distribution to determine the probability of an agent updating and communicating at each time step. That is, this probability determines what time steps are in the sets $K^i$ and $R_j^i$. We consider a range of probabilities in the interval $[0.1, 1]$, where a probability 1 of updating

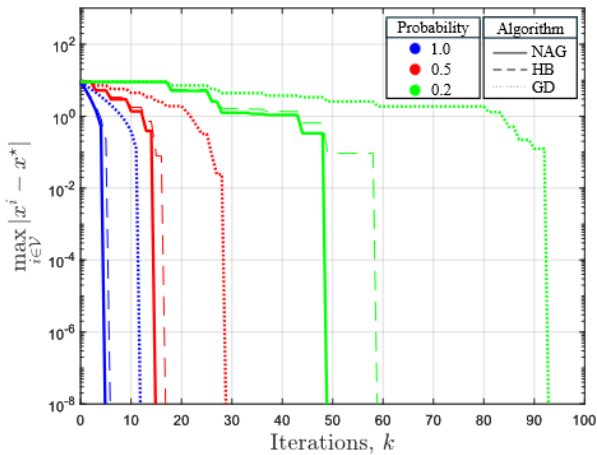

Fig. 1. The evolution of the worst-performing agent in the totally asynchronous NAG algorithm (solid lines), heavy ball algorithm (dashed lines), and gradient descent (dotted lines). In all cases, the NAG algorithm that we have developed converges faster than both the heavy ball algorithm and a gradient descent algorithm that are run with the same computation and communication times.

TABLE I

OPERATION COMPLEXITY OF NAG, HEAVY BALL (HB), AND GRADIENT DESCENT (GD).

| Probability | Iterations | | | Percent Reduction | |
|---|---|---|---|---|---|
| | NAG | HB | GD | HB | GD |
| 1.0 | 5 | 6 | 12 | 17% | 58% |
| 0.9 | 7 | 7 | 15 | 0% | 53% |
| 0.8 | 10 | 11 | 17 | 9% | 41% |
| 0.7 | 8 | 10 | 20 | 20% | 60% |
| 0.6 | 10 | 11 | 24 | 9% | 58% |
| 0.5 | 15 | 17 | 29 | 12% | 48% |
| 0.4 | 15 | 19 | 34 | 21% | 56% |
| 0.3 | 24 | 29 | 58 | 17% | 59% |
| 0.2 | 49 | 59 | 93 | 17% | 47% |
| 0.1 | 123 | 170 | 314 | 28% | 61% |

and communicating leads to the double-step synchronous algorithm of Section III-B.

Figure 1 displays the distance to the optimum for three probabilities in $[0.1, 1]$, and Table I provides the operation complexity results for a wider range of probabilities. Both demonstrate the faster convergence rate of Algorithm 1 compared to the totally asynchronous heavy ball algorithm in [15] and totally asynchronous gradient descent. In simulations, it was seen that heavy ball converged in, at best, the same time as NAG, though often slower, as Figure 1 and Table I show. Gradient descent was seen to converge substantially slower than both NAG and heavy ball. In addition, Figure 1 shows that NAG outpaces gradient descent and heavy ball by a wider margin as the probabilities of updating and communicating decrease. We observe a maximum reduction in convergence time of $28\%$ between NAG and heavy ball, and a maximum reduction in convergence time of $61\%$ between NAG and gradient descent. These results indicate the superior performance of NAG relative to existing algorithms in the totally asynchronous context.

## VII. CONCLUSION

This paper presented what is, to the best of our knowledge, the first totally asynchronous implementation of Nesterov's accelerated gradient algorithm for optimization. We showed that this algorithm converges linearly in the number of agents' computations and communications when counted in a certain sequence, and simulations showed that it converges faster than comparable heavy ball and gradient descent algorithms. Future work will identify additional forms of accelerated algorithms that converge under total asynchrony and explore the implementation of these techniques in teams of mobile robots whose communications are subject to jamming.

REFERENCES

[1] J. C. Derenick and J. R. Spletzer, "Convex optimization strategies for coordinating large-scale robot formations," *IEEE Transactions on Robotics*, vol. 23, no. 6, pp. 1252–1259, 2007.
[2] R. Deits and R. Tedrake, "Computing large convex regions of obstacle-free space through semidefinite programming," in *Algorithmic Foundations of Robotics XI: Selected Contributions of the Eleventh International Workshop on the Algorithmic Foundations of Robotics*. Springer, 2015, pp. 109–124.
[3] J. Alonso-Mora, S. Baker, and D. Rus, "Multi-robot navigation in formation via sequential convex programming," in *2015 IEEE/RSJ International Conference on Intelligent Robots and Systems (IROS)*. IEEE, 2015, pp. 4634–4641.
[4] S. Bubeck *et al.*, "Convex optimization: Algorithms and complexity," *Foundations and Trends® in Machine Learning*, vol. 8, no. 3-4, pp. 231–357, 2015.
[5] H. Jiang, Y. T. Lee, Z. Song, and S. C.-w. Wong, "An improved cutting plane method for convex optimization, convex-concave games, and its applications," in *Proceedings of the 52nd Annual ACM SIGACT Symposium on Theory of Computing*, 2020, pp. 944–953.
[6] J. Mattingley and S. Boyd, "Real-time convex optimization in signal processing," *IEEE Signal processing magazine*, vol. 27, no. 3, pp. 50–61, 2010.
[7] D. P. Palomar and Y. C. Eldar, *Convex optimization in signal processing and communications*. Cambridge university press, 2010.
[8] P. Tseng, D. P. Bertsekas, and J. N. Tsitsiklis, "Partially asynchronous, parallel algorithms for network flow and other problems," *SIAM Journal on Control and Optimization*, vol. 28, no. 3, pp. 678–710, 1990.
[9] D. Bertsekas and J. Tsitsiklis, *Parallel and distributed computation: numerical methods*. Athena Scientific, 2015.
[10] K. R. Hendrickson and M. T. Hale, "Towards totally asynchronous primal-dual convex optimization in blocks," in *2020 59th IEEE Conference on Decision and Control (CDC)*. IEEE, 2020, pp. 3663–3668.
[11] M. Ubl and M. Hale, "Totally asynchronous large-scale quadratic programming: Regularization, convergence rates, and parameter selection," *IEEE Transactions on Control of Network Systems*, vol. 8, no. 3, pp. 1465–1476, 2021.
[12] K. R. Hendrickson and M. T. Hale, "Totally asynchronous primal-dual convex optimization in blocks," *IEEE Transactions on Control of Network Systems*, vol. 10, no. 1, pp. 454–466, 2023.
[13] S. Hochhaus and M. T. Hale, "Asynchronous distributed optimization with heterogeneous regularizations and normalizations," in *IEEE Conference on Decision and Control (CDC)*, 2018, pp. 4232–4237.
[14] F. Mansoori and E. Wei, "A fast distributed asynchronous newton-based optimization algorithm," *IEEE Transactions on Automatic Control*, vol. 65, no. 7, pp. 2769–2784, 2019.
[15] D. M. Hustig-Schultz, K. Hendrickson, M. Hale, and R. G. Sanfelice, "A totally asynchronous block-based heavy ball algorithm for convex optimization," in *2023 American Control Conference (ACC)*. IEEE, 2023, pp. 873–878.
[16] D. M. Hustig-Schultz and R. G. Sanfelice, "Uniting nesterov and heavy ball methods for uniform global asymptotic stability of the set of minimizers," *Automatica*, vol. 160, p. 111389, 2024.
[17] R. T. Rockafellar and R. J.-B. Wets, *Variational analysis*. Springer Science & Business Media, 2009, vol. 317.

## APPENDIX

### A. Proof of Lemma 1

A necessary and sufficient condition for a point $x^\star \in \mathcal{X}$ to be the minimizer of $f$ over $\mathcal{X}$ is to have $\langle x - x^\star, -\nabla f(x) \rangle \leq 0$ for all $x \in \mathcal{X}$ by Proposition 3.1 in [9, pg. 210]. By Proposition 5.7 in [9, pg. 275], a vector $x^\star$ is optimal if and only if $\langle x_i^i - x_i^\star, -\gamma \nabla_i f(x^\star) \rangle \leq 0$ for every $x_i^i \in \mathcal{X}_i$ and $i \in \mathcal{V}$.

We can change the second argument of the inner product using

$$-\gamma \nabla_i f(x^\star) = x_i^\star - \gamma \nabla_i f\big(x^\star + \lambda(x^\star - x^\star)\big) + \lambda(x_i^\star - x_i^\star) - x_i^\star,$$

where $\gamma, \lambda \in \mathbb{R}_+$. Then we have

$$\langle x_i^i - x_i^\star, x_i^\star - \gamma \nabla_i f\big(x^\star + \lambda(x^\star - x^\star)\big) + \lambda(x_i^\star - x_i^\star) - x_i^\star \rangle \leq 0, \quad (16)$$

for all $x_i^i \in \mathcal{X}_i$. For a function $f$ that is twice continuously differentiable and convex, and a set $\mathcal{X}$ that is nonempty, compact, and convex, and some $v \in \mathbb{R}^n$, a vector $u \in \mathcal{X}$ satisfies $u = \Pi_\mathcal{X}[v]$ if and only if $\langle y - u, v - u \rangle \leq 0$ for all $y \in \mathcal{X}$ by Proposition 3.2 of [9, pg. 211]. Thus, by (16),

$$x_i^\star = \Pi_{\mathcal{X}_i}\left[x_i^\star - \gamma \nabla_i f\big(x^\star + \lambda(x^\star - x^\star)\big) + \lambda(x_i^\star - x_i^\star)\right].$$

Therefore, for all $i \in \mathcal{V}$, the point $(x^\star, x^\star)$ is a fixed point for the update law in (4) since $x_i^\star = \tilde{u}_x^i(x^\star, x^\star)$.

Consider the update law in (5), namely $y_i^i(l+1) = \tilde{u}_y^i\big(x^i(l), y^i(l)\big) = x_i^i(l)$. Here, the point $(x^\star, x^\star)$ is a fixed point of this equation for all $i \in \mathcal{V}$, since $x_i^\star = \tilde{u}_y^i(x^\star, x^\star)$.

Therefore, $z_i^\star \in \mathcal{Z}_i$ is a fixed point of the update law $\tilde{u}^i$ in (6) for all $i \in \mathcal{V}$ because $z_i^\star = \tilde{u}^i(z^\star)$.

### B. Proof of Theorem 1

We have established with Lemma 1 that $z_i^\star = \tilde{u}^i(z^\star)$ for all $i \in \mathcal{V}$, i.e., we have both $x_i^\star = \tilde{u}_x^i(z^\star)$ and $x_i^\star = \tilde{u}_y^i(z^\star)$ for all $i \in \mathcal{V}$. Next, using the definition of the $\infty$-norm and (4), we have

$$\|x(l+1) - x^\star\|_\infty = \max_{i \in \mathcal{V}} \left| x_i^i(l+1) - x_i^\star \right|$$
$$= \max_{i \in \mathcal{V}} \Big| \Pi_{\mathcal{X}_i}\left[ x_i^i(l) - \gamma \nabla_i f\big(x^i(l) + \lambda(x^i(l) - y^i(l))\big) \right.$$
$$\left. + \lambda(x_i^i(l) - y_i^i(l)) \right] - \Pi_{\mathcal{X}_i}\left[ z_i^\star - \gamma \nabla_i f\big(z^\star + \lambda(z^\star - z^\star)\big) \right.$$
$$\left. + \lambda(z_i^\star - z_i^\star) \right] \Big|,$$

where we have used the aforementioned fixed point property.

The projection operator is nonexpansive with respect to $|\cdot|$ by Proposition 3.2 in [9, pg. 211], i.e., $|\Pi_{\mathcal{X}_i}[v_2] - \Pi_{\mathcal{X}_i}[v_1]| \leq |v_2 - v_1|$ for all $v_1, v_2 \in \mathbb{R}$. Applying that property, we find

$$\|x(l+1) - x^\star\|_\infty \leq \max_{i \in \mathcal{V}} \Big| x_i^i(l) - \gamma \nabla_i f\big(x^i(l)$$
$$+ \lambda(x^i(l) - y^i(l))\big) + \lambda(x_i^i(l) - y_i^i(l))$$
$$- \Big(z_i^\star - \gamma \nabla_i f\big(z^\star + \lambda(z^\star - z^\star)\big) + \lambda(z_i^\star - z_i^\star)\Big)\Big|.$$

We can then use (7) to write

$$\|x(l+1) - x^\star\|_\infty \leq \max_{i \in \mathcal{V}} \left| \hat{u}^i_{\text{true}}(z^{\text{true}}(l)) - \hat{u}^i_{\text{true}}(z^\star) \right|. \quad (17)$$

From the Mean Value Theorem [9, pg. 639], we know that for a continuously differentiable function $g : \mathbb{R}^n \to \mathbb{R}$, there exists some $\rho \in [0, 1]$ such that

$$g(v) - g(w) = \langle \nabla g(\rho w + (1-\rho)v), v - w \rangle,$$

where $v, w \in \mathbb{R}^n$ are vectors.

For each $i \in \mathcal{V}$, we will apply the Mean Value Theorem to the map $\hat{u}^i_{\text{true}}$ defined in (7), for which we define $a^i = (a_1^i, a_2^i)$, $c^i = (c_1^i, c_2^i)$, and $b^i = (b_1^i, b_2^i)$ where the vectors $a_1^i, a_2^i, c_1^i, c_2^i, b_1^i, b_2^i \in \mathbb{R}^n$. Then the $j^{\text{th}}$ component of $b_1^i$ satisfies $b_{1,j}^i \in [a_{1,j}^i, c_{1,j}^i]$ and the $j^{\text{th}}$ component of $b_2^i$ satisfies $b_{2,j}^i \in [a_{2,j}^i, c_{2,j}^i]$. Then, applying the Mean Value Theorem to $\hat{u}^i_{\text{true}}$ and expanding, we find

$$\hat{u}^i_{\text{true}}(c^i) - \hat{u}^i_{\text{true}}(a^i) = \qquad (18)$$
$$= \sum_{j=1}^n \frac{\partial \hat{u}^i_{\text{true}}(b^i)}{\partial x_j}(c_{1,j}^i - a_{1,j}^i) + \sum_{j=1}^n \frac{\partial \hat{u}^i_{\text{true}}(b^i)}{\partial y_j}(c_{2,j}^i - a_{2,j}^i).$$

For $i \neq j$, the required partial derivatives of $\hat{u}^i_{\text{true}}(b^i)$ are

$$\frac{\partial \hat{u}^i_{\text{true}}(b^i)}{\partial x_i} = 1 - \gamma(1+\lambda)\nabla_i^2 f\big(b_1^i + \lambda(b_1^i - b_2^i)\big) + \lambda$$
$$\frac{\partial \hat{u}^i_{\text{true}}(b^i)}{\partial x_j} = -\gamma(1+\lambda)\nabla_j \nabla_i f\big(b_1^i + \lambda(b_1^i - b_2^i)\big)$$
$$\frac{\partial \hat{u}^i_{\text{true}}(b^i)}{\partial y_i} = \gamma \lambda \nabla_i^2 f\big(b_1^i + \lambda(b_1^i - b_2^i)\big) - \lambda$$
$$\frac{\partial \hat{u}^i_{\text{true}}(b^i)}{\partial y_j} = \gamma \lambda \nabla_j \nabla_i f\big(b_1^i + \lambda(b_1^i - b_2^i)\big).$$

Substituting these expressions back into (18) results in

$$\hat{u}^i_{\text{true}}(c^i) - \hat{u}^i_{\text{true}}(a^i)$$
$$= (1 - \gamma(1+\lambda)\nabla_i^2 f\big(b_1^i + \lambda(b_1^i - b_2^i)\big) + \lambda)(c_{1,i}^i - a_{1,i}^i)$$
$$+ \sum_{\substack{j=1 \\ j \neq i}}^n (-\gamma(1+\lambda)\nabla_j \nabla_i f\big(b_1^i + \lambda(b_1^i - b_2^i)\big))(c_{1,j}^i - a_{1,j}^i)$$
$$+ (\gamma \lambda \nabla_i^2 f\big(b_1^i + \lambda(b_1^i - b_2^i)\big) - \lambda)(c_{2,i}^i - a_{2,i}^i)$$
$$+ \sum_{\substack{j=1 \\ j \neq i}}^n (\gamma \lambda \nabla_j \nabla_i f\big(b_1^i + \lambda(b_1^i - b_2^i)\big))(c_{2,j}^i - a_{2,j}^i).$$

Taking the absolute value of both sides and applying the triangle inequality to the right-hand side yields

$$\left| \hat{u}^i_{\text{true}}(c^i) - \hat{u}^i_{\text{true}}(a^i) \right|$$
$$\leq \left| 1 - \gamma(1+\lambda)\nabla_i^2 f(b_1^i + \lambda(b_1^i - b_2^i)) + \lambda \right| \left| c_{1,i}^i - a_{1,i}^i \right|$$
$$+ \gamma \sum_{\substack{j=1 \\ j \neq i}}^n \left| -(1+\lambda)\nabla_j \nabla_i f(b_1^i + \lambda(b_1^i - b_2^i)) \right| \left| c_{1,j}^i - a_{1,j}^i \right|$$
$$+ \left| \gamma \lambda \nabla_i^2 f(b_1^i + \lambda(b_1^i - b_2^i)) - \lambda \right| \left| c_{2,i}^i - a_{2,i}^i \right|$$

$$+ \gamma \sum_{\substack{j=1\\j\neq i}}^{n} \left|\lambda \nabla_j \nabla_i f(b_1^i + \lambda(b_1^i - b_2^i))\right| \left|c_{2,j}^i - a_{2,j}^i\right|.$$

Replacing $\nabla_i^2 f$ with $H_{ii}$ and $\nabla_j \nabla_i$ with $H_{ij}$, where $H = \nabla^2 f$ is the Hessian of $f$, we reach

$$\left|\hat{u}_{\text{true}}^i(c^i) - \hat{u}_{\text{true}}^i(a^i)\right| \quad (19)$$
$$\leq \left|1 - \gamma(1+\lambda)H_{ii}(b_1^i + \lambda(b_1^i - b_2^i)) + \lambda\right| \left|c_{1,i}^i - a_{1,i}^i\right|$$
$$+ \gamma \sum_{\substack{j=1\\j\neq i}}^{n} \left|-(1+\lambda)H_{ij}(b_1^i + \lambda(b_1^i - b_2^i))\right| \left|c_{1,j}^i - a_{1,j}^i\right|$$
$$+ \left|\gamma\lambda H_{ii}(b_1^i + \lambda(b_1^i - b_2^i)) - \lambda\right| \left|c_{2,i}^i - a_{2,i}^i\right|$$
$$+ \gamma \sum_{\substack{j=1\\j\neq i}}^{n} \left|\lambda H_{ij}(b_1^i + \lambda(b_1^i - b_2^i))\right| \left|c_{2,j}^i - a_{2,j}^i\right|.$$

We now define $\gamma$ to ensure that the terms that mulitply the $\left|c_{1,i}^i - a_{1,i}^i\right|$ term and the $\left|c_{2,i}^i - a_{2,i}^i\right|$ term in (19) are equal to their absolute values. To do so, we take $\gamma > 0$ and

$$\gamma < \frac{1}{\max_{i\in\mathcal{V}} \max_{\eta\in\mathcal{X}} |H_{ii}(\eta)|},$$

where we leverage the fact that every diagonal entry of the Hessian is positive, i.e., $H_{ii}(\eta) > 0$ for all $\eta \in \mathcal{X}$, due to Assumption 3. Then, using

$$\gamma H_{ii}(b_1^i + \lambda(b_1^i - b_2^i)) < \frac{H_{ii}(b_1^i + \lambda(b_1^i - b_2^i))}{\max_{i\in\mathcal{V}} \max_{\eta\in\mathcal{X}} |H_{ii}(\eta)|} \leq 1,$$

we arrive at

$$1 - \gamma(1+\lambda)H_{ii}(b_1^i + \lambda(b_1^i - b_2^i)) + \lambda \geq 0.$$

A similar procedure can be established for the term that multiplies $\left|c_{2,i}^i - a_{2,i}^i\right|$, though in this case we establish that the negated version of that term is nonnegative, i.e.,

$$-\lambda\gamma H_{ii}(b_1^i + \lambda(b_1^i - b_2^i)) + \lambda \geq 0.$$

Then, (19) becomes

$$\left|\hat{u}_{\text{true}}^i(c^i) - \hat{u}_{\text{true}}^i(a^i)\right| \quad (20)$$
$$\leq (1 - \gamma(1+\lambda)H_{ii}(b_1^i + \lambda(b_1^i - b_2^i)) + \lambda) \left|c_{1,i}^i - a_{1,i}^i\right|$$
$$+ \gamma(1+\lambda) \sum_{\substack{j=1\\j\neq i}}^{n} \left|H_{ij}(b_1^i + \lambda(b_1^i - b_2^i))\right| \left|c_{1,j}^i - a_{1,j}^i\right|$$
$$+ (-\gamma\lambda H_{ii}(b_1^i + \lambda(b_1^i - b_2^i)) + \lambda) \left|c_{2,i}^i - a_{2,i}^i\right|$$
$$+ \gamma\lambda \sum_{\substack{j=1\\j\neq i}}^{n} \left|H_{ij}(b_1^i + \lambda(b_1^i - b_2^i))\right| \left|c_{2,j}^i - a_{2,j}^i\right|.$$

Since $\left\|c_1^i - a_1^i\right\|_\infty = \max_{j\in\mathcal{V}} \left|c_{1,j}^i - a_{1,j}^i\right|$ and $\left\|c_2^i - a_2^i\right\|_\infty = \max_{j\in\mathcal{V}} \left|c_{2,j}^i - a_{2,j}^i\right|$, we can rewrite (20) as

$$\left|\hat{u}_{\text{true}}^i(c^i) - \hat{u}_{\text{true}}^i(a^i)\right| \quad (21)$$
$$\leq \Bigg(1 - \gamma(1+\lambda)H_{ii}(b_1^i + \lambda(b_1^i - b_2^i)) + \lambda$$
$$+ \gamma(1+\lambda) \sum_{\substack{j=1\\j\neq i}}^{n} \left|H_{ij}(b_1^i + \lambda(b_1^i - b_2^i))\right|\Bigg) \left\|c_1^i - a_1^i\right\|_\infty$$
$$+ \Bigg(-\gamma\lambda H_{ii}(b_1^i + \lambda(b_1^i - b_2^i)) + \lambda$$
$$+ \gamma\lambda \sum_{\substack{j=1\\j\neq i}}^{n} \left|H_{ij}(b_1^i + \lambda(b_1^i - b_2^i))\right|\Bigg) \left\|c_2^i - a_2^i\right\|_\infty.$$

Further, due to Assumption 3, it is true that $H_{ii}(\eta) - \sum_{\substack{j=1\\j\neq i}} |H_{ij}(\eta)| \geq \mu$ for $\mu > 0$ and all $\eta \in \mathcal{X}$ and all $i \in \mathcal{V}$. Applying this inequality, we upper bound (21) as

$$\left|\hat{u}_{\text{true}}^i(c^i) - \hat{u}_{\text{true}}^i(a^i)\right| \leq (1+\lambda-\gamma\mu(1+\lambda)) \left\|c_1^i - a_1^i\right\|_\infty$$
$$+ \lambda(1-\gamma\mu) \left\|c_2^i - a_2^i\right\|_\infty. \quad (22)$$

Now we set $c^i = z^{\text{true}}(l)$ and $a^i = z^\star$, which also gives $c_1^i = x^{\text{true}}(l)$, $c_2^i = y^{\text{true}}(l)$, and $a_1^i = a_2^i = z^\star$. With these substitutions, we use (22) in (17) to find

$$\|x(l+1) - x^\star\|_\infty \leq (1+\lambda-\gamma\mu(1+\lambda)) \left\|x^{\text{true}}(l) - x^\star\right\|_\infty$$
$$+ \lambda(1-\gamma\mu) \left\|y^{\text{true}}(l) - x^\star\right\|_\infty \quad (23)$$

for all $l \in \mathbb{N}$. This bound holds for all $l \in \mathbb{N}$, and therefore

$$\|x(l) - x^\star\|_\infty \leq (1+\lambda-\gamma\mu(1+\lambda)) \left\|x^{\text{true}}(l-1) - x^\star\right\|_\infty$$
$$+ \lambda(1-\gamma\mu) \left\|y^{\text{true}}(l-1) - x^\star\right\|_\infty. \quad (24)$$

We apply (24) in (23) and use $y^{\text{true}}(l) = x^{\text{true}}(l-1)$ to find

$$\|x(l+1) - x^\star\|_\infty \leq$$
$$((1+\lambda-\gamma\mu(1+\lambda))^2 + \lambda(1-\gamma\mu)) \left\|x^{\text{true}}(l-1) - x^\star\right\|_\infty$$
$$+ \lambda(1+\lambda-\gamma\mu(1+\lambda))(1-\gamma\mu) \left\|y^{\text{true}}(l-1) - x^\star\right\|_\infty. \quad (25)$$

Next, using $y(l+1) = x(l)$, we use (24) see that

$$\|y(l+1) - x^\star\|_\infty \leq \|x(l) - x^\star\|_\infty$$
$$\leq (1+\lambda-\gamma\mu(1+\lambda)) \left\|x^{\text{true}}(l-1) - x^\star\right\|_\infty$$
$$+ \lambda(1-\gamma\mu) \left\|y^{\text{true}}(l-1) - x^\star\right\|_\infty \quad (26)$$

for all $l \in \mathbb{N}$.

Because $z^{\text{true}}(l-1) = (x^{\text{true}}(l-1), y^{\text{true}}(l-1))$, the definition of the infinity-norm implies that

$$\left\|x^{\text{true}}(l-1) - x^\star\right\|_\infty \leq \left\|z^{\text{true}}(l-1) - z^\star\right\|_\infty$$
$$\left\|y^{\text{true}}(l-1) - x^\star\right\|_\infty \leq \left\|z^{\text{true}}(l-1) - z^\star\right\|_\infty.$$

Therefore, (25) and (26) can be bounded as

$$\|x(l+1) - x^\star\|_\infty \leq \Big((1+\lambda-\gamma\mu(1+\lambda))^2 + \lambda(1-\gamma\mu)$$
$$+ \lambda(1-\gamma\mu)(1+\lambda-\gamma\mu(1+\lambda))\Big) \left\|z^{\text{true}}(l-1) - z^\star\right\|_\infty$$

and

$$\|y(l+1) - x^\star\|_\infty \leq (1 - \gamma\mu + 2\lambda(1 - \gamma\mu)) \|z^{\text{true}}(l-1) - z^\star\|_\infty.$$

We define

$$\alpha_1 = (1 + \lambda - \gamma\mu(1+\lambda))^2 + \lambda(1-\gamma\mu) \\ + \lambda(1-\gamma\mu)(1+\lambda-\gamma\mu(1+\lambda)) \quad (27)$$

and

$$\alpha_2 = 1 - \gamma\mu + 2\lambda(1-\gamma\mu) \quad (28)$$

and reach

$$\|z(l+1) - z^\star\|_\infty \leq \alpha \|z^{\text{true}}(l-1) - z^\star\|_\infty,$$

for all $l \in \mathbb{N}$ and $i \in \mathcal{V}$, where $\alpha = \max\{\alpha_1, \alpha_2\} \geq 0$. In the synchronous setting, $z(l) = z^{\text{true}}(l)$, due to synchronous communication after each computation is performed. Therefore, we have

$$\|z(l+1) - z^\star\|_\infty \leq \alpha \|z(l-1) - z^\star\|_\infty.$$

Lastly, we will define the hyperparameter $\lambda$ to ensure that $\alpha_1, \alpha_2 < 1$. Due to the fact that $\gamma > 0$, $\mu > 0$, and $H_{ii}(\eta) \geq \mu$ for all $\eta \in \mathcal{X}$ and $i \in \mathcal{V}$, we take

$$0 < \lambda < \frac{\gamma\mu}{2(1-\gamma\mu)}.$$

Substituting this into (27) and (28) gives $\alpha \in [0, 1)$.

### C. Proof of Lemma 4

From the SCC in Lemma 3, we see that $z \in \mathcal{Z}(k)$ implies that $h(z) \in \mathcal{Z}(k+1)$ for all $k \in \mathbb{N}$. Recall that $h$ maps the global state of the network to the successive true state by applying the local update $u^i$ for all $i \in \mathcal{V}$. In the totally asynchronous setting, it can take arbitrarily long for all agents to update their local decision variables as stated in Assumption 4. However, for $z^i(k) \in \mathcal{Z}(k)$, when agent $i$ has computed $z_i^i(k+1)$, it is then true that $z_i^i(k+1) \in \mathcal{Z}_i(k+1)$. Since $\mathcal{Z}_i(k+1) \subset \mathcal{Z}_i(k)$ by the LLC in Lemma 3, we still have $z_i^i(k+1) \in \mathcal{Z}_i(k)$ and thus agent $i$ still has $z^i(k+1) \in \mathcal{Z}(k)$, even if no other entries of $z^i$ have changed onboard agent $i$.

Therefore, $\mathcal{Z}_i(k)$ is forward invariant under agent $i$'s computations for all $k \in \mathbb{N}$ and for all agents $i \in \mathcal{V}$. Since this statement applies for computations with each agent's local decision variables, it also holds for all communications of those variables within in the network. In other words, it holds that $\mathcal{Z}_j(k)$ is forward invariant for all $j \in \mathcal{V}^i$ under both computations and communications. Therefore, because each set $\mathcal{Z}_i(k)$ is forward invariant, by the Box Condition, we conclude that the full set $Z(k) = \mathcal{Z}_1(k) \times \mathcal{Z}_2(k) \times \cdots \times \mathcal{Z}_n(k)$ is forward invariant.